\documentclass[12pt]{article}
\usepackage[cp1251]{inputenc}
\usepackage[english]{babel}
\usepackage{amsfonts}

\newtheorem{theor}{Theorem}
\newtheorem{prop}{Proposition}
\newtheorem{lem}{Lemma}

\begin{document}

\title{Cones of weighted quasi-metrics, weighted quasi-hypermetrics and of 
oriented cuts}
\author{M.Deza \and V.Grishukhin \and E.Deza}
\date{}
\maketitle

\begin{abstract}
We show that the cone of weighted $n$-point quasi-metrics $WQMet_n$, the cone of weighted quasi-hypermetrics $WHyp_n$ and the cone of oriented cuts $OCut_n$ are projections along an extreme ray of the metric cone $Met_{n+1}$, of the hypermetric cone $Hyp_{n+1}$ and of the cut cone $Cut_{n+1}$, respectively. This projection is such that if one knows all faces of an original cone then one knows all faces of the projected cone. \end{abstract}

\section{Introduction}
Oriented (or directed) distances are encountered very often, for 
example, these are one-way transport routes, rivers with quick flow and so on.

The notions of directed distances, quasi-metrics and oriented cuts are 
generalizations of the notions of distances, metrics and cuts, respectively (see, for example, \cite{DL}),  which are 
central objects in Graph Theory and Combinatorial Optimization.

Quasi-metrics are used in Semantics of Computations (see, for example, \cite{Se}) and in computational geometry (see, for example, \cite{AAC}).  Oriented distances have been used already by Hausdorff in 1914, see \cite{Ha}.

In \cite{CMM}, authors give an example of directed metric derived from a metric as follows. Let $d$ be a metric on a set $V\cup\{0\}$, where 0 is a distinguished point. Then a quasi-metric $q$ on the set $V$ is given as
\[q_{ij}=d_{ij}+d_{i0}-d_{j0}.\]

This quasi-metric belongs to a special important subclass of quasi-metrics, 
namely, to a class of {\em weighted quasi-metrics}. We show (cf. also Lemma 1 (ii) in \cite{DDV})  
that {\em any} weighted quasi-metric is obtained by a slight generalization 
of this method.

All semi-metrics on a set of cardinality $n$ form a {\em metric cone} $Met_n$.
 There are two important sub-cones of $Met_n$, namely, the cone $Hyp_n$ of 
{\em hypermetrics}, and the cone $Cut_n$ of {\em $\ell_1$-metrics}. 
These three cones form the following nested family $Cut_n\subset Hyp_n\subset Met_n$, see \cite{DL}.

 We introduce a space $Q_n$, called a {\em space of weighted 
quasi-metrics} and define in it  a cone $WQMet_n$. Elements of this cone 
satisfy triangle and non-negativity inequalities. Among extreme rays of 
the cone $WQMet_n$ there are rays spanned by {\em ocut vectors}, i.e., incidence vectors of oriented cuts.

We define in the space $Q_n$ a cone $OCut_n$ as the cone hull of ocut vectors. Elements of the cone $OCut_n$ are weighted quasi-$\ell$-metrics.

Let semi-metrics in the cone $Met_{n+1}$ be defined on the set $V\cup\{0\}$. 
The {\em cut cone} $Cut_{n+1}$ (or the {\em cone of  $\ell_1$-metrics} on 
this set is a cone hull of cut semi-metrics $\delta(S)$ for all 
$S\subset V\cup\{0\}$. The cut semi-metrics $\delta(S)$ are extreme rays of all the three cones $Met_{n+1}$, $Hyp_{n+1}$ and $Cut_{n+1}$. In particular, $\delta(\{0\})=\delta(V)$ is an extreme ray of these three cones.

In this paper, it is shown that the cones $WQMet_n$ and $OCut_n$ are 
projections of the corresponding cones $Met_{n+1}$ and 
$Cut_{n+1}$ along the extreme ray $\delta(V)$. We define a cone 
$WQHyp_n$ of {\em weighted quasi-hypermetrics} as projection along $\delta(V)$ of the cone $Hyp_{n+1}$. So, we obtain a nested family $OCut_n\subset WQHyp_n \subset WQMet_n$.

The cones of weighted quasi-metrics, oriented cuts  and other related  
generalizations of metrics are studied in \cite{DD} and \cite{DDV}. 
The  polytope of oriented cuts was considered in \cite{AM}.

\section{Spaces ${\mathbb R}^E$ and ${\mathbb R}^{E^{\cal O}}$}
\label{spR}
Let $V$ be a set of cardinality $|V|=n$. Let $E$ and $E^{\cal O}$ be sets of all unordered $(ij)$ and ordered $ij$ pairs of elements $i,j\in V$. Consider two Euclidean spaces ${\mathbb R}^E$ and ${\mathbb R}^{E^{\cal O}}$of vectors $d\in{\mathbb R}^E$ and $g\in{\mathbb R}^{E^{\cal O}}$ with coordinates $d_{(ij)}$ and $g_{ij}$, where $(ij)\in E$ and $ij\in E^{\cal O}$, respectively. Obviously, dimensions of the spaces ${\mathbb R}^E$ and ${\mathbb R}_s^{E^{\cal O}}$ are $|E|=\frac{n(n-1)}{2}$ and $|E^{\cal O}|=n(n-1)$, respectively.

Denote by $(d,t)=\sum_{(ij)\in E}d_{(ij)}t_{(ij)}$ scalar product of vectors 
$d,t\in{\mathbb R}^E$. Similarly, $(f,g)=\sum_{ij\in E^{\cal O}}f_{ij}g_{ij}$ 
 is the scalar product of vectors $f,g\in{\mathbb R}^{E^{\cal O}}$.

Let $\{e_{(ij)}:(ij)\in E\}$ and $\{e_{ij}:ij\in {E^{\cal O}}\}$ be orthonormal bases of ${\mathbb R}^E$ and ${\mathbb R}^{E^{\cal O}}$, respectively. Then, for $f\in{\mathbb R}^E$ and $q\in{\mathbb R}^{E^{\cal O}}$, we have
\[(e_{(ij)},f)=f_{(ij)}\mbox{  and  }(e_{ij},q)=q_{ij}. \]

For $f\in{\mathbb R}^{E^{\cal O}}$, define $f^*\in{\mathbb R}^{E^{\cal O}}$ as follows
\[f^*_{ij}=f_{ji}\mbox{  for all  }ij\in E^{\cal O}. \]
Each vector $g\in{\mathbb R}^{E^{\cal O}}$ can be decompose into {\em symmetric} $g^s$ and {\em antisymmetric} $g^a$ parts as follows:
\[g^s=\frac{1}{2}(g+g^*), \mbox{  }g^a=\frac{1}{2}(g-g^*), \mbox{  }g=g^s+g^a. \]
Call a vector $g$ {\em symmetric} if $g^*=g$, and {\em antisymmetric} if $g^*=-g$. Let ${\mathbb R}_s^{E^{\cal O}}$ and ${\mathbb R}_a^{E^{\cal O}}$ be subspaces of the corresponding vectors.
Note that the spaces ${\mathbb R}_s^{E^{\cal O}}$ and ${\mathbb R}_a^{E^{\cal O}}$ are mutually orthogonal. In fact, for $p\in{\mathbb R}_s^{E^{\cal O}}$ and $f\in{\mathbb R}_a^{E^{\cal O}}$, we have
\[(p,f)=\sum_{ij\in E^{\cal O}}p_{ij}f_{ij}=\sum_{(ij)\in E}(p_{ij}f_{ij}
+p_{ji}f_{ji})=\sum_{(ij)\in E}(p_{ij}f_{ij}-p_{ij}f_{ij})=0. \]
Hence
\[{\mathbb R}^{E^{\cal O}}={\mathbb R}_s^{E^{\cal O}}\oplus{\mathbb R}_a^{E^{\cal O}}, \]
where $\oplus$ is the direct sum.

Obviously, there is an isomorphism $\varphi$ between the spaces ${\mathbb R}^E$ and ${\mathbb R}_s^{E^{\cal O}}$. Let $d\in{\mathbb R}^E$ have coordinates $d_{(ij)}$. Then
\[d^{\cal O}=\varphi(d)\in{\mathbb R}_s^{E^{\cal O}}\mbox{ ,  such that  } d^{\cal O}_{ij}=d^{\cal O}_{ji}=d_{(ij)}. \]
In particular,
\[\varphi(e_{(ij)})=e_{ij}+e_{ji}. \]

The map $\varphi$ is invertible. In fact, for 
$q\in{\mathbb R}_s^{E^{\cal O}}$, we have $\varphi^{-1}(q)=d\in{\mathbb R}^E$, 
such that $d_{(ij)}=q_{ij}=q_{ji}$.
The isomorphism $\varphi$ will be useful in what follows.

\section{Space of weights $Q_n^w$}
\label{Wsp}
One can consider the sets $E$ and $E^{\cal O}$ as sets of edges $(ij)$ and arcs $ij$ of an unordered and ordered complete graphs $K_n$ and $K_n^{\cal O}$ on the vertex set $V$, respectively. The graph $K_n^{\cal O}$ has two arcs $ij$ and $ji$ between each pair of vertices $i,j\in V$.

It is convenient to consider vectors $g\in{\mathbb R}^{E^{\cal O}}$ as functions on the set of arcs $E^{\cal O}$ of the graph $K_n^{\cal O}$. So, the decomposition ${\mathbb R}^{E^{\cal O}}={\mathbb R}_s^{E^{\cal O}}\oplus{\mathbb R}_a^{E^{\cal O}}$ is a decomposition of the space of all functions on arcs in $E^{\cal O}$ onto the spaces of symmetric and antisymmetric functions.

Besides, there is an important direct decomposition of the 
space ${\mathbb R}_a^{E^{\cal O}}$ of antisymmetric functions 
into two subspaces. In the Theory of Electric Networks, these spaces are called 
 spaces of {\em tensions} and {\em flows} (see also \cite{Aig}).

The tension space relates to {\em potentials} (or {\em weights}) $w_i$ given on vertices $i\in V$ of the graph $K_n^{\cal O}$. The corresponding antisymmetric function $g^w$ is determined as
\[g^w_{ij}=w_i-w_j.\]
It is called {\em tension} on the arc $ij$. Obviously, 
$g^w_{ji}=w_j-w_i=-g^w_{ij}$. Denote by $Q_n^w$ the subspace of 
${\mathbb R}^{E^{\cal O}}$ generated by all tensions on arcs 
$ij\in E^{\cal O}$. We call $Q_n^w$  a {\em space of weights}.

Each tension function $g^w$ is represented as weighted sum of 
elementary {\em potential} functions $q(k)$, for $k\in V$, as follows:
\[g^w=\sum_{k\in V}w_kq(k), \]
where
\begin{equation}
\label{qek}
q(k)=\sum_{j\in V-\{k\}}(e_{kj}-e_{jk}),\mbox{  for all  }k\in V,
\end{equation}
are basic functions that generate the space of weights $Q_n^w$.
Hence, the values of the basic functions $q(k)$ on arcs are as follows:
\begin{equation}
\label{qk}
q_{ij}(k)=\left\{
\begin{array}{rl}
1, &\mbox{ if }i=k\\
-1, &\mbox{ if }j=k\\
0,& \mbox{ otherwise.}
\end{array}\right.
\end{equation}
We obtain
 \[g^w_{ij}=\sum_{k\in V}w_kq_{ij}(k)=w_i-w_j.\]
It is easy to verify that
\[q^2(k)=(q(k),q(k))=2(n-1), \mbox{  }(q(k),q(l))=-2\mbox{  for all } k,l\in V, k\not=l, \mbox{  }\sum_{k\in V}q(k)=0. \]
Hence, there are only $n-1$ independent functions $q(k)$ that generate the space $Q_n^w$.

The weighted quasi-metrics lie in the space ${\mathbb R}_s^{E^{\cal O}}\oplus Q_n^w$ that we denote as $Q_n$. Direct complements of  $Q_n^w$ in ${\mathbb R}_a^{E^{\cal O}}$ and $Q_n$ in ${\mathbb R}^{E^{\cal O}}$ is a space $Q_n^c$ of {\em circuits} (or {\em flows}).

\section{Space of circuits $Q_n^c$}
\label{Csp}
The {\em space of circuits} (or {\em space of flows}) is generated by 
characteristic vectors of oriented circuits in the graph 
$K_n^{\cal O}$. Arcs of $K_n^{\cal O}$ are ordered pairs 
$ij$ of vertices $i,j\in V$. The arc $ij$ is oriented from the vertex $i$ to the vertex $j$. Recall that $K_n^{\cal O}$ has both the arcs $ij$ and $ji$ for each pair of vertices $i,j\in V$.

Let $G_s\subset K_n$ be a subgraph with a set of edges $E(G_s)\subset E$. 
We relate to the graph $G_s$ a directed graph 
$G\subset K_n^{\cal O}$ with the arc set $E^{\cal O}(G)\subset E^{\cal O}$ as follows. An arc $ij$ belongs to $G$, i.e., $ij\in E^{\cal O}(G)$, 
if and only if $(ij)=(ji)\in E(G)$.  This definition implies that in 
this case,  the arc $ji$ belongs to $G$ also, i.e., $ji\in E^{\cal O}(G)$.

Let $C_s$ be a circuit in the graph $K_n$. The circuit $C_s$ is determined 
by a sequence of distinct vertices $i_k\in V$, where $1\le k\le p$, and 
$p$ is the length of $C_s$. The edges of $C_s$ are unordered 
pairs $(i_k,i_{k+1})$, where indices are taken modulo $p$. By above 
definition, an {\em oriented bicircuit} $C$ of the graph $K_n^{\cal O}$ 
relates  to the circuit $C_s$. Arcs of $C$ are ordered pairs $i_ki_{k+1}$ and $i_{k+1}i_k$, where indices are taken modulo $p$. Take an orientation of $C$. Denote by $-C$ the {\em opposite} circuit with opposite orientation. Denote an arc of $C$ {\em direct} or {\em opposite} if its direction coincides with or is opposite to the given orientation of $C$, respectively. Let $C^+$ and $C^-$ be subcircuits of $C$ consisting of direct and opposite arcs, respectively.

The following vector $f^C$ is the characteristic vector of the bicircuit $C$:
\[f^C_{ij}=\left\{
\begin{array}{rl}
1, &\mbox{ if } ij\in C^+,   \\
-1. &\mbox{ if } ij\in C^-, \\
0,& \mbox{ otherwise.}
\end{array}\right. \]

Note that $f^{-C}=(f^C)^*=-f^C$, and $f^C\in{\mathbb R}_a^{E^{\cal O}}$.

Denote by $Q_n^c$ the space linearly generated by circuit vectors $f^C$ for all bicircuits $C$ of the graph $K_n^{\cal O}$. It is well known that characteristic vectors of {\em fundamental} circuits form a basis of $Q_n^c$. Fundamental circuits are defined as follows.

Let $T$ be a spanning tree of the graph $K_n$. Since $T$ is spanning, its vertex set $V(T)$ is the set of all vertices of $K_n$, i.e., $V(T)=V$. Let $E(T)\subset E$ be the set of edges of $T$. Then any edge $e=(ij)\not\in E(T)$ closes a unique path in $T$ between vertices $i$ and $j$ into a circuit $C_s^e$. This circuit $C_s^e$ is called {\em fundamental}. Call corresponding oriented bicircuit $C^e$ also {\em fundamental}.

There are $|E-E(T)|=\frac{n(n-1)}{2}-(n-1)$ fundamental circuits. Hence
 \[{\rm dim}Q_n^c=\frac{n(n-1)}{2}-(n-1), \mbox{  and  }{\rm dim}Q_n+{\rm dim}Q_n^c=n(n-1)={\rm dim}{\mathbb R}^{E^{\cal O}}.\]
This implies that $Q_n^c$ is an orthogonal complement of $Q_n^w$ in ${\mathbb R}_a^{\cal O}$ and $Q_n$ in ${\mathbb R}^{E^{\cal O}}$, i.e.
\[{\mathbb R}_a^{E^{\cal O}}=Q_n^w\oplus Q_n^c\mbox{  and  }{\mathbb R}^{E^{\cal O}}=Q_n\oplus Q_n^c={\mathbb R}^{E_s^{\cal O}}\oplus Q_n^w\oplus Q_n^c.\]

\section{Cut and ocut vector set-functions}
\label{Cut}
The space $Q_n$ is generated also by vectors of oriented cuts, which we define in this section.

Each subset $S\subset V$ determines cuts of the graphs $K_n$ and $K_n^{\cal O}$ that are subsets of edges and arcs of these graphs.

A $cut(S)\subset E$ is a subset of edges $(ij)$ of $K_n$ such that $(ij)\in cut(S)$ if and only if $|\{i,j\}\cap S|=1$.

A $cut^{\cal O}(S)\subset E^{\cal O}$ is a subset of arcs $ij$ of $K_n^{\cal O}$ such that $ij\in cut^{\cal O}(S)$ if and only if $|\{i,j\}\cap S|=1$. So, if $ij\in cut^{\cal O}(S)$, then $ji\in cut^{\cal O}(S)$ also.

An {\em oriented cut} is a subset $ocut(S)\subset E^{\cal O}$  of arcs $ij$ of $K_n^{\cal O}$ such that $ij\in ocut(S)$ if and only if $i\in S$ and $j\not\in S$.

We relate to these three types of cuts characteristic vectors $\delta(S)\in{\mathbb R}^E$, $\delta^{\cal O}(S)\in{\mathbb R }_s ^{E^{\cal O}}$, $q(S)\in{\mathbb R}_a^{E^{\cal O}}$ and $c(S)\in{\mathbb R}^{E^{\cal O}}$ as follows.

For $cut(S)$, we set
\[\delta(S)=\sum_{i\in S,j\in{\overline S}}e_{(ij)},\mbox{  such that  }
\delta_{(ij)}(S)=\left\{
\begin{array}{rl}
1, &\mbox{ if }|\{i,j\}\cap S|=1\\
0,& \mbox{ otherwise,}
\end{array}\right.\]
where ${\overline S}=V-S$.
For $cut^{\cal O}(S)$, we set
\[\delta^{\cal O}(S)=\varphi(\delta(S))=\sum_{i\in S,j\in{\overline S}}(e_{ij}+e_{ji})
\mbox{  and  }q(S)=\sum_{i\in S,j\in{\overline S}}(e_{ij}-e_{ji}). \]
Hence,
\[\delta^{\cal O}_{ij}(S)=\left\{
\begin{array}{rl}
1, &\mbox{ if }|\{i,j\}\cap S|=1\\
0,& \mbox{ otherwise.}
\end{array}\right.
\mbox{  and  }
q_{ij}(S)=\left\{
\begin{array}{rl}
1, &\mbox{ if }i\in S, j\not\in S\\
-1, &\mbox{ if }j\in S, i\not\in S\\
0,& \mbox{ otherwise.}
\end{array}\right. \]
Note that, for one-element sets $S=\{k\}$, the function $q(\{k\})$ is $q(k)$ of section 2. It is easy to see that
\[(\delta^{\cal O}(S),q(T))=0\mbox{  for any  }S,T\subset V.\]

For the oriented cut $ocut(S)$, we set
\[c(S)=\sum_{i\in S,j\in{\overline S}}e_{ij}. \]
Hence,
\[c_{ij}(S)=\left\{
\begin{array}{rl}
1, &\mbox{ if }i\in S, j\not\in S\\
0,& \mbox{ otherwise.}
\end{array}\right. \]
Obviously, it holds $c(\emptyset)=c(V)={\bf 0}$, where ${\bf 0}\in{\mathbb R}^{E^{\cal O}}$ is a vector whose all coordinates are equal zero.
We have
\begin{equation}
\label{dec}
c^*(S)=c({\overline S}), \mbox{   }c(S)+c({\overline S})=\delta^{\cal O}(S), \mbox{  }c(S)-c({\overline S})=q(S)\mbox{  and  }c(S)=\frac{1}{2}(\delta^{\cal O}(S)+q(S)).
\end{equation}
Besides, we have
\[c^s(S)=\frac{1}{2}\delta^{\cal O}(S), \mbox{   }c^a(S)=\frac{1}{2}q(S). \]

Recall that a set-function $f(S)$ on all $S\subset V$, is called {\em submodular} if, for any $S,T\subset V$, the following 
{\em submodular inequality} holds: 
\[f(S)+f(T)-(f(S\cap T)+f(S\cup T))\ge 0. \]
It is well known that the vector set-function $\delta\in{\mathbb R}^E$ is submodular (see, for example, \cite{Aig}). The above isomorphism $\varphi$   of the spaces ${\mathbb R}^E$ and ${\mathbb R}_s^{E^{\cal O}}$ implies that the vector set-function $\delta^{\cal O}=\varphi(\delta)\in{\mathbb R}_s^{E^{\cal O}}$ is submodular also.

A set-function $f(S)$ is called {\em modular} if, for any $S,T\subset V$, the above submodular inequality holds as equality. This equality is called {\em modular equality}. It is well known (and can be easily verified) that antisymmetric vector set-function $f^a(S)$ is modular for  any oriented graph $G$. Hence, our antisymmetric vector set-function $q(S)\in{\mathbb R}_a^{E^{\cal O}}$ for the oriented complete graph $K_n^{\cal O}$ is modular also.

Note that the set of all submodular set-functions on a set $V$ forms a 
cone in the space ${\mathbb R}^{2^V}$. Therefore, the last equality in (\ref{dec}) implies that the vector set-function $c(S)\in{\mathbb R}^{E^{\cal O}}$ is submodular.

The modularity of the antisymmetric vector set-function $q(S)$ is important for what follows. It is well-known (see, for example, \cite{Bir}) (and it can be easily verified using modular equality) that a modular set-function $m(S)$ is completely determined by its values on the empty set and on all one-element sets. Hence, a modular set-function $m(S)$ has the following form
\[ m(S)=m_0+\sum_{i\in S}m_i, \]
where $m_0=m(\emptyset)$ and $m_i=m(\{i\})-m(\emptyset)$. For brevity, we 
set $f(\{i\})=f(i)$ for any set-function $f(S)$. Since $q(\emptyset)=q(V)=0$, we have
\begin{equation}
\label{qS}
q(S)=\sum_{k\in S}q(k),\mbox{   }S\subset V, \mbox{  and  }q(V)=\sum_{k\in V}q(k)=0.
\end{equation}
Using equations (\ref{dec}) and (\ref{qS}), we obtain
\begin{equation}
\label{cS}
c(S)=\frac{1}{2}(\delta^{\cal O}(S)+\sum_{k\in S}q(k)).
\end{equation}

Now we show that ocut vectors $c(S)$ for all $S\subset V$ linearly generate the space  $Q_n\subseteq{\mathbb R}^{E^{\cal O}}$. The space generated by $c(S)$ consists of the following vectors
\[c=\sum_{S\subset V}\alpha_Sc(S), \mbox{  where  }\alpha_S\in{\mathbb R}. \]

Recall that $c(S)=\frac{1}{2}(\delta^{\cal O}(S)+q(S))$. Hence, we have
\[c=\frac{1}{2}\sum_{S\subset V}\alpha_S(\delta^{\cal O}(S)+q(S))=
\frac{1}{2}\sum_{S\subset V}\alpha_S\delta^{\cal O}(S)+ \frac{1}{2}\sum_{S\subset V}\alpha_Sq(S)=\frac{1}{2}(d^{\cal O}+q),\]
where $d^{\cal O}=\varphi(d)$ for $d=\sum_{S\subset V}\alpha_S\delta(S)$. 
For a vector $q$, we have
\[q=\sum_{S\subset V}\alpha_Sq(S)=\sum_{S\subset V}\alpha_S\sum_{k\in S}q(k)=\sum_{k\in V}w_kq(k),\mbox{  where  }
w_k=\sum_{V\supset S\ni k}\alpha_S. \]

Since $q_{ij}=\sum_{k\in V}w_kq_{ij}(k)=w_i-w_j$, we have

\begin{equation}
\label{cw}
 c_{ij}=\frac{1}{2}(d^{\cal O}_{ij}+w_i-w_j).
 \end{equation}

It is well-known (see, for example, \cite{DL}) that the cut vectors 
$\delta(S)\in{\mathbb R}^E$ for all $S\subset V$ linearly generate 
the full space ${\mathbb R}^E$. Hence,  the vectors 
$\delta^{\cal O}(S)\in{\mathbb R}_s^{E^{\cal O}}$, for all $S\subset V$, linearly generate the full space ${\mathbb R}_s^{E^{\cal O}}$.

According to (\ref{cS}), antisymmetric parts of ocut vectors $c(S)$ generate the space $Q_n^w$. This implies that the space $Q_n={\mathbb R}_s^{E^{\cal O}}\oplus Q_n^w$ is generated by $c(S)$ for all $S\subset V$.

\section{Properties of the space $Q_n$}
\label{prop}
Let $x\in Q_n$ and let $f^C$ be the characteristic vector of a bicircuit $C$. Since $f^C$ is orthogonal to $Q_n$, we have $(x,f^C)=\sum_{ij\in C}f^C_{ij}x_{ij}=0$. This equality implies that each point $x\in Q_n$ satisfies the following equalities
\[\sum_{ij\in C^+}x_{ij}=\sum_{ij\in C^-}x_{ij} \]
for any bicircuit $C$.

Let $K_{1,n-1}\subset K_n$ be a spanning star of $K_n$ consisting of all $n-1$ edges incident to a vertex of $K_n$. Let this vertex be 1. Each edge of $K_n-K_{1,n-1}$ has the form $(ij)$, where $i\not=1\not=j$. The edge $(ij)$ closes a fundamental triangle with edges $(1i),(1j),(ij)$. The corresponding bitriangle $T(1ij)$ generates the equality
\[x_{1i}+x_{ij}+x_{j1}=x_{i1}+x_{1j}+x_{ji}. \]
These equalities are the case $k=3$ of
{\em k-cyclic symmetry} considered in \cite{DD}. They were  derived by another way in 
\cite{AM}. They correspond to fundamental bi-triangles $T(1ij)$, for all $i,j\in V-\{1\}$, and are all $\frac{n(n-1)}{2}-(n-1)$ independent equalities determining the space, where the $Q_n$ lies.

Above coordinates $x_{ij}$ of a vector $x\in Q_n$ are given in the orthonormal basis $\{e_{ij}:ij\in E^{\cal O}\}$. But, for what follows, it is more convenient to consider vectors $q\in Q_n$ in another basis. Recall that ${\mathbb R}_s^{E^{\cal O}}=\varphi({\mathbb R}^E)$. Let, for $(ij)\in E$, $\varphi(e_{(ij)})=e_{ij}+e_{ji}\in{\mathbb R}_s^{E^{\cal O}}$ be basic vectors of the subspace ${\mathbb R}_s^{E^{\cal O}}\subset Q_n$. Let $q(i)\in Q_n^w$, $i\in V$, be basic vectors of the space $Q_n^w\subset Q_n$. Then, for $q\in Q_n$, we set
\[q=q^s+q^a,\mbox{  where  }q^s=\sum_{(ij)\in E}q_{(ij)}\varphi(e_{(ij)}), \mbox{  }q^a=\sum_{i\in V}w_iq(i). \]

Now, we obtain an important expression for the scalar product $(g,q)$ of vectors $g,q\in Q_n$. Recall that $(\varphi(e_{(ij)}),q(k))= ((e_{ij}+e_{ji}),q(k))=0$ for all $(ij)\in E$ and all $k\in V$. Hence
$(g^s,q^a)=(g^a,q^s)=0$, and we have
\[(g,q)=(g^s,q^s)+(g^a,q^a).\]
Besides, we have
\[((e_{ij}+e_{ji}),(e_{kl}+e_{lk}))=0\mbox{  if  }(ij)\not=(kl), \mbox{  }(e_{ij}+e_{ji})^2=2, \]
and (see Section~\ref{Wsp})
\[(q(i),q(j))=-2\mbox{  if  }i\not=j,\mbox{  }(q(i))^2=2(n-1). \]

Let $v_i$, $i\in V$, be weights of the vector $g$. Then we have
\[(g,q)=2\sum_{(ij)\in E}g_{(ij)}q_{(ij)}+2(n-1)\sum_{i\in V}v_iw_i-2\sum_{i\not=j\in V}v_iw_j. \]
For the last sum, we have
\[\sum_{i\not=j\in V}v_iw_j=(\sum_{i\in V}v_i)(\sum_{i\in V}w_i)- \sum_{i\in V}v_iw_i. \]
Since weights are defined up to an additive scalar, we can choose weights $v_i$ such that $\sum_{i\in V}v_i=0$. Then the last sum in the product $(g,q)$ is equal to $-\sum_{i\in V}v_iw_i$. Finally we obtain that the sum of antisymmetric parts is equal to $2n\sum_{i\in V}v_iw_i$. So, for the product of two vectors $g,q\in Q_n$ we have the following expression
\[(g,q)=(g^s,q^s)+(g^a,q^a)=2(\sum_{(ij)\in E}g_{(ij)}q_{(ij)}+n\sum_{i\in V}v_iw_i)\mbox{  if  }\sum_{i\in V}v_i=0\mbox{  or  }\sum_{i\in V}w_i=0. \]
In what follows, we consider inequalities $(g,q)\ge 0$. We can delete the multiple 2, and rewrite such inequality as follows
\begin{equation}
\label{prod}
\sum_{(ij)\in E}g_{(ij)}q_{(ij)}+n\sum_{i\in V}v_iw_i\ge 0,
\end{equation}
where $\sum_{i\in V}v_i=0$.

Below we consider some cones in the space $Q_n$. Since the space $Q_n$ is orthogonal to the space of circuits $Q_n^c$, each facet vector of a cone in $Q_n$ is defined up to a vector of the space $Q_n^c$. Of course each vector $g'\in{\mathbb R}^{E^{\cal O}}$ can be decomposed as $g'=g+g^c$, where $g\in Q_n$ and $g^c\in Q_n^c$. Call the vector $g\in Q_n$ {\em canonical representative} of the vector $g'$. Usually we will use canonical facet vectors. But sometimes not canonical representatives of a facet vector are useful.

Cones $Con$ that will be considered are invariant under the 
operation $q\to q^*$, defined in Section~\ref{spR}. In other words, $Con^*=Con$. This operation changes signs of weights:
\[q_{ij}=q_{(ij)}+w_i-w_j\to q_{(ij)}+w_j-w_i=q_{(ij)}-w_i+w_j.\]

Let $(g,q)\ge 0$ be an inequality determining a facet $F$ of a cone $Con\subset 
Q_n$. Since $Con=Con^*$, the cone $Con$ has, together with the facet $F$, also a facet $F^*$. The facet $F^*$ is determined by the inequality $(g^*,q)\ge 0$.

\section{Projections of cones $Con_{n+1}$}
Recall that $Q_n={\mathbb R}_s^{E^{\cal O}}\oplus Q_n^w$, ${\mathbb R}_s^{E^{\cal O}}=\varphi({\mathbb R}^E)$ and dim$Q_n=\frac{n(n+1)}{2}-1$. Let $0\not\in V$ be an additional point. Then the set of unordered pairs $(ij)$ for $i,j\in V\cup\{0\}$ is $E\cup E_0$, where $E_0=\{(0i):i\in V\}$. Obviously, ${\mathbb R}^{E\cup E_0}={\mathbb R}^E\oplus{\mathbb R}^{E_0}$ and dim${\mathbb R}^{E\cup E_0}=\frac{n(n+1)}{2}$.

The space ${\mathbb R}^{E\cup E_0}$ contains the following three important cones: the cone $Met_{n+1}$ of semi-metrics, the cone $Hyp_{n+1}$ of hyper-semi-metrics and the cone $Cut_{n+1}$ of $\ell_1$-semi-metrics, all on the set $V\cup\{0\}$. Denote by $Con_{n+1}$ any of these cones.

Recall that a semi-metric $d=\{d_{(ij)}\}$ is called {\em metric} if $d_{(ij)}\not=0$ for all $(ij)\in E$. For brevity sake, in what follows, we call elements of the cones $Con_{n+1}$ simply metrics (or hypermetrics, $\ell_1$-metrics), assuming that they can be semi-metrics.

Note that if $d\in Con_{n+1}$ is a metric on the set $V\cup\{0\}$, then a restriction $d^V$ of $d$ on the set $V$ is a point of the cone $Con_n=Con_{n+1}\cap{\mathbb R}^E$ of metrics on the set $V$. In other words, we can suppose that $Con_n\subset Con_{n+1}$.

The cones $Met_{n+1}$, $Hyp_{n+1}$ and $Cut_{n+1}$ contain the cut vectors  $\delta(S)$ that span extreme rays for all $S\subset V\cup\{0\}$. Denote by $l_0$ the extreme ray spanned by the cut vector $\delta(V)= \delta(\{0\})$. Consider a projection $\pi({\mathbb R}^{E\cup E_0})$ of the space ${\mathbb R}^{E\cup E_0}$ along the ray $l_0$ onto a subspace of ${\mathbb R}^{E\cup E_0}$ that is orthogonal to $\delta(V)$. This projection is such that $\pi({\mathbb R}^E)={\mathbb R}^E$ and $\pi({\mathbb R}^{E\cup E_0})={\mathbb R}^E\oplus\pi({\mathbb R}^{E_0})$.

Note that $\delta(V)\in{\mathbb R}^{E_0}$, since, by Section~\ref{Cut}, $\delta(V)=\sum_{i\in V}e_{(0i)}$. For simplicity sake, set
\[e_0=\delta(\{0\})=\delta(V)=\sum_{i\in V}e_{(0i)}. \]
Recall that the vector $e_0$ spans the extreme ray $l_0$. Obviously, the space ${\mathbb R}^E$ is orthogonal to $l_0$, and therefore, $\pi({\mathbb R}^E)={\mathbb R}^E$.

Let $x\in {\mathbb R}^E$. We decompose this point as follows
\[x=x^V+x^0, \]
where $x^V=\sum_{(ij)\in E}x_{(ij)}e_{(ij)}\in{\mathbb R}^E$ and $x^0=\sum_{i\in V}x_{(0i)}e_{(0i)}\in{\mathbb R}^{E_0}$.
We define a map $\pi$ as follows:
\[\pi(e_{(ij)})=e_{(ij)}\mbox{  for  }(ij)\in E,\mbox{  and } \pi(e_{(0i)})=e_{(0i)}-\frac{1}{n}e_0\mbox{  for  }i\in V. \]

So, we have
\begin{equation}
\label{px}
\pi(x)=\pi(x^V)+\pi(x^0)=\sum_{(ij)\in E}x_{(ij)}e_{(ij)}+\sum_{i\in V}x_{(0i)} (e_{(0i)}-\frac{1}{n}e_0).
\end{equation}

Note that the projection $\pi$ transforms the positive orthant of the space ${\mathbb R}^{E_0}$ onto the whole space $\pi({\mathbb R}^{E_0})$.

Now we describe how faces of a cone in the space ${\mathbb R}^{E\cup E_0}$ are projected along one of its extreme rays.

Let $l$ be an extreme ray and $F$ be a face of a cone in ${\mathbb R}^{E\cup E_0}$. Let $\pi$ be the projection along $l$. Let {\rm dim}$F$ be dimension of the face $F$. Then the following equality holds
\begin{equation}
\label{Fl}
{\rm dim}\pi(F)={\rm dim}F-{\rm dim}(F\cap l).
\end{equation}
Let $g\in{\mathbb R}^{E\cup E_0}$ be a facet vector of a facet $G$, and $e$ be a vector spanning the line $l$. Then dim$(G\cap l)=1$ if $(g,e)=0$, and dim$(G\cap l)=0$ if $(g,e)\not=0$.

\begin{theor}
\label{fac}
Let $G$ be a face of the cone $\pi(Con_{n+1})$. Then $G=\pi(F)$, where $F$ is a face of $Con_{n+1}$ such that there is a facet of $Con_{n+1}$, containing $F$ and the extreme ray $l_0$ spanned by $e_0=\delta(V)$.

In particular, $G$ is a facet of $\pi(Con_{n+1})$ if and only if  $G=\pi(F)$, where $F$ is a facet of $Con_{n+1}$ containing the extreme ray $l_0$. Similarly, $l'$ is an extreme ray of $\pi(Con_{n+1})$ if and only if $l'=\pi(l)$, where $l$ is an extreme ray of $Con_{n+1}$ lying in a facet of $Con_{n+1}$ that contains $l_0$.
\end{theor}
{\bf Proof}. Let $\cal F$ be a set of all facets of the cone $Con_{n+1}$. Then $\cup_{f\in{\cal F}}\pi(F)$ is a covering of the projection $\pi(Con_{n+1})$. By (\ref{Fl}), in this covering, if $l_0\subset F\in{\cal F}$, then $\pi(F)$ is a facet of $\pi(Con_{n+1})$. If $l_0\not\subset F$, then there is a one-to-one correspondence between points of $F$ and $\pi(F)$. Hence, dim$\pi(F)=n$, and $\pi(F)$ cannot be a facet of $\pi(Con_{n+1})$, since $\pi(F)$ fills an $n$-dimensional part of the cone $\pi(Con_{n+1})$.

If $F'$ is a face of $Con_{n+1}$, then $\pi(F')$ is a face of the above covering. If $F'$ belongs only to facets $F\in{\cal F}$ such that $l_0\not\subset F$, then $\pi(F')$ lies inside of $\pi(Con_{n+1})$. In this case, it is not a face of $\pi(Con_{n+1})$. This implies that $\pi(F')$ is a face of $\pi(Con_{n+1})$ if and only if $F'\subset F$, where $F$ is a facet of $Con_{n+1}$ such that $l_0\subset F$. Suppose that dimension of $F'$ is $n-1$, and $l_0\not\subset F'$. Then dim$\pi(F')=n-1$. If $F'$ is contained in a facet $F$ of $Con_{n+1}$ such that $l_0\subset F$, then $\pi(F')=\pi(F)$. Hence, $\pi(F')$ is a facet of the cone $\pi(Con_{n+1})$ that coincides with the facet $\pi(F)$.

Now, the assertions of Theorem about facets and extreme rays of $\pi(Con_{n+1})$ follow. \hfill $\Box$

\vspace{2mm}
Theorem~\ref{fac} describes all faces of the cone $\pi(Con_{n+1})$ if one knows all faces of the cone $Con_{n+1}$.

Recall that we consider $Con_n=Con_{n+1}\cap{\mathbb R}^E$ as a sub-cone of $Con_{n+1}$, and therefore, $\pi(Con_n)\subset\pi(Con_{n+1})$. Since $\pi({\mathbb R}^E)={\mathbb R}^E$, we have $\pi(Con_n)=Con_n$. Let $(f,x)\ge $ be a facet-defining inequality of a facet $F$ of the cone $Con_{n+1}$. Since $Con_{n+1}\subset{\mathbb R}^E\oplus{\mathbb R}^{E_0}$, we represent vectors $f,x\in{\mathbb R}^{E\cup E_0}$ as $f=f^V+f^0,x=x^V+x^0$, where $f^V,x^V\in{\mathbb R}^E$ and $f^0,x^0\in{\mathbb R}^{E_0}$. Hence, the above facet-defining inequality can be rewritten as
\[(f,x)=(f^V,x^V)+(f^0,x^0)\ge 0. \]

It turns out that $Con_{n+1}$ has always a facet $F$ with its facet vector $f=f^V+f^0$ such that $f^0=0$. Since $f^V$ is orthogonal to ${\mathbb R}^{E_0}$, the hyperplane $(f^V,x)=(f^V,x^V)=0$ supporting the facet $F$ contains the whole space ${\mathbb R}^{E_0}$. The equality $(f^V,x^V)=0$ defines a facet $F^V=F\cap{\mathbb R}^E$ of the cone $Con_n$.

{\bf Definition}. A facet $F$ of the cone $Con_{n+1}$ with a facet vector $f=f^V+f^0$ is called {\em zero-lifting} of a facet $F^V$ of $Con_n$ if $f^0=0$ and $F\cap{\mathbb R}^E=F^V$.

Similarly, a facet $\pi(F)$ of the cone $\pi(Con_{n+1})$ with a facet vector $f$ is called zero-lifting of $F^V$ if $f=f^V$ and $\pi(F)\cap{\mathbb R}^E=F^V$.

\vspace{2mm}
It is well-known (see, for example, \cite{DL}) that each facet $F^V$ with facet vector $f^V$ of the cone $Con_n$ can be zero-lifted up to a facet $F$ of $Con_{n+1}$ with the same facet vector $f^V$.

\begin{prop}
\label{piF}
Let a facet $F$ of $Con_{n+1}$ be a zero-lifting of a facet $F^V$ of $Con_n$. Then $\pi(F)$ is a facet of $\pi(Con_{n+1})$ that is also zero-lifting of $F^V$.
\end{prop}
{\bf Proof}. Recall that the hyperplane $\{x\in{\mathbb R}^{E\cup E_0}:(f^V,x)=0\}$ supporting the facet $F$ contains the whole space ${\mathbb R}^{E_0}$. Hence, the facet $F$ contains the extreme ray $l_0$ spanned by the vector $e_0\in{\mathbb R}^{E_0}$. By Theorem~\ref{fac}, $\pi(F)$ is a facet of $\pi(Con_{n+1})$. The facet vector of $\pi(F)$ can be written as $f=f^V+f'$, where $f^V\in{\mathbb R}^E$ and $f'\in\pi({\mathbb R}^{E_0})$. Since the hyperplane supporting the facet $\pi(F)$ is given by the equality $(f^V,x)=0$ for $x\in\pi({\mathbb R}^{E\cup E_0})$, we have $f'=0$. Besides, obviously, $\pi(F)\cap{\mathbb R}^E=F^V$. Hence, $\pi(F)$ is zero-lifting of $F^V$. \hfill $\Box$

\section{Cones $\psi(Con_{n+1})$}
\label{psi}
Note that basic vectors of the space ${\mathbb R}^{E\cup E_0}$ are $e_{(ij)}$ for $(ij)\in E$ and $e_{(0i)}$ for $(0i)\in E_0$.
Since $\pi(e_0)=\sum_{i\in V}\pi(e_{(0i)})=0$, we have dim$\pi({\mathbb R}^{E_0})=n-1=$dim$Q_n^w$. Note that $\pi({\mathbb R}^E)={\mathbb R}^E$. Hence, there is a one-to-one bijection $\chi$ between the spaces $\pi({\mathbb R}^{E\cup E_0})$ and $Q_n$.

We define this bijection $\chi:\pi({\mathbb R}^{E\cup E_0})\to Q_n$ as follows
\[\chi({\mathbb R}^E)=\varphi({\mathbb R}^E)={\mathbb R}_s^{E^{\cal O}}, \mbox{  and  }\chi(\pi({\mathbb R}^{E_0}))=Q_n^w, \]
where
\[\chi(e_{(ij)})=\varphi(e_{(ij)})=e_{ij}+e_{ji}, \mbox{  and  } \chi(\pi(e_{(0i)}))=\chi(e_{(0i)}-\frac{1}{n}e_0)=q(i), \]
where $q(i)$ is defined in (\ref{qek}).

Note that $(e_{ij}+e_{ji})^2=2=2e_{(ij)}^2$ and
\[(q(i),q(j))=-2=2n((e_{(0i)}-\frac{1}{n}e_0),(e_{(0j)}-\frac{1}{n}e_0)), \mbox{   }q^2(i)=2(n-1)=2n(e_{(0i)}-\frac{1}{n}e_0)^2. \]
Roughly speaking, the map $\chi$ is a homothety that extends vectors $e_{(0i)}-\frac{1}{n}e_0$ up to vectors $q(i)$ by the multiple $\sqrt{2n}$.

Setting $\psi=\chi\circ\pi$, we obtain a map $\psi:{\mathbb R}^{E\cup E_0}\to Q_n$ such that
\begin{equation}
\label{pse}
\psi(e_{(ij)})=e_{ij}+e_{ji}\mbox{  for  }(ij)\in E, \mbox{  }\psi(e_{(0i)})=q(i)\mbox{  for  }i\in V.
\end{equation}

Now we show how a point $x=x^V+x^0\in{\mathbb R}^{E\cup E_0}$ is transformed into a point $q=\psi(x)=\chi(\pi(x))\in Q_n$. We have  $\pi(x)=x^V+\pi(x^0)$, where, according to (\ref{px}), $x^V=\sum_{(ij)\in E}x_{(ij)}e_{(ij)}\in\pi({\mathbb R}^E)={\mathbb R}^E$ and $\pi(x^0)=\sum_{i\in V}x_{(0i)}(e_{(0i)}-\frac{1}{n}e_0)\in\pi({\mathbb R}^{E_0})$. Obviously, $\chi(x^V+\pi(x^0))= \chi(x^V)+\chi(\pi(x^0))$, and
\[\psi(x^V)=\chi(x^V)=\sum_{(ij)\in E}x_{(ij)}(e_{ij}+e_{ji})=\varphi(x^V)=q^s\mbox{  and  } \chi(\pi(x^0))=\sum_{i\in V}x_{(0i)}q(i)=q^a. \]
Recall that $q^s=\sum_{(ij)\in E}q_{(ij)}(e_{ij}+e_{ji})$ and $q^a=\sum_{i\in V}w_iq(i)$. Hence
\begin{equation}
\label{coo}
q_{(ij)}=x_{(ij)},\mbox{  }(ij)\in E,\mbox{  and  }w_i=x_{(0i)}. \end{equation}

Let $f\in{\mathbb R}^{E\cup E_0}$ be a facet vector of a facet $F$ of the cone $Con_{n+1}$, $f=f^V+f^0=\sum_{(ij)\in E}f_{(ij)}e_{(ij)}+\sum_{i\in V}f_{(0i)}e_{(0i)}$.

Let $(f,x)\ge 0$ be the inequality determining the facet $F$.
The inequality $(f,x)\ge 0$ takes on the set $V\cup\{0\}$ the following form
\[(f,x)=\sum_{(ij)\in E}f_{(ij)}x_{(ij)}+\sum_{i\in V}f_{(0i)}x_{(0i)}\ge 0. \]
Since $x_{(ij)}=q_{(ij)}$, $x_{(0i)}=w_i$, we can rewrite this inequality as follows
\begin{equation}
\label{gF}
(f,q)=(f^V,q^s)+(f^0,q^a)\equiv\sum_{(ij)\in E}f_{(ij)}q_{(ij)}+\sum_{i\in V}f_{(0i)}w_i\ge 0.
\end{equation}

Comparing the inequality (\ref{gF}) with (\ref{prod}), we see that a canonical form of the facet vector $f$ is $f=f^s+f^a$, where
\begin{equation}
\label{fv}
f^s_{(ij)}=f_{(ij)},\mbox{  for  }(ij)\in E, \mbox{  } f^a_{ij}=v_i-v_j\mbox{  where  }v_i=\frac{1}{n}f_{(0i)},\mbox{  }i\in V.
\end{equation}

\begin{theor}
\label{main}
Let $F$ be a facet of the cone $Con_{n+1}$. Then $\psi(F)$ is a facet of the cone $\psi(Con_{n+1})$ if and only if the facet $F$ contains the extreme ray $l_0$ spanned by the vector $e_0$.

Let $l\not=l_0$ be be an extreme ray of $Con_{n+1}$. Then $\psi(l)$ is an extreme ray of $\psi(Con_{n+1})$ if and only if the ray $l$ belongs to a facet containing the extreme ray $l_0$.
\end{theor}
{\bf Proof}. By Theorem~\ref{fac}, the projection $\pi$ transforms the facet $F$ of $Con_{n+1}$ into a facet of $\pi(Con_{n+1})$ if and only if $l_0\subset F$. By the same Theorem, the projection $\pi(l)$ is an extreme ray of $\pi(Con_{n+1})$ if and only if $l$ belongs to a facet containing the extreme ray $l_0$.

Recall that the map $\chi$ is a bijection between the spaces ${\mathbb R}^{E\cup E_0}$ and $Q_n$. This implies the assertion of this Theorem for the map $\psi=\chi\circ\pi$. \hfill $\Box$

\vspace{2mm}
By Theorem~\ref{main}, the map $\psi$ transforms the facet $F$ in a facet of the cone $\psi(Con_{n+1})$ only if $F$ contains the extreme ray $l_0$, 
i.e., only if the equality $(f,e_0)=0$ holds. Hence, the facet vector $f$ should satisfy the equality $\sum_{i\in V}f_{(0i)}=0$.

The inequalities (\ref{gF}) give all facet-defining inequalities of the cone $\psi(Con_{n+1})$ from known facet-defining inequalities of the cone $Con_{n+1}$.

A proof of Proposition~\ref{fso} below will be given later for each of the cones $Met_{n+1}$, $Hyp_{n+1}$ and $Cut_{n+1}$ separately.
\begin{prop}
\label{fso}
Let $F$ be a facet of $Con_{n+1}$ with facet vector $f=f^V+f^0$ such that $(f^0,e_0)=0$. Then $Con_{n+1}$ has also a facet $F^*$ with facet vector $f^*=f^V-f^0$.
\end{prop}

\vspace{2mm}
Proposition~\ref{fso} implies the following important fact.
\begin{prop}
\label{iso}
For $q=q^s+q^a\in\psi(Con_{n+1})$, the map $q=q^s+q^a\to q^*=q^s-q^a$ preserves the cone $\psi(Con_{n+1})$, i.e.
\[(\psi(Con_{n+1}))^*=\psi(Con_{n+1}). \]
\end{prop}
{\bf Proof}. Let $F$ be a facet of $Con_{n+1}$ with facet vector $f$. By Proposition~\ref{fso}, if $\psi(F)$ is a facet of $\psi(Con_{n+1})$, then $F^*$ is a facet of $Con_{n+1}$ with facet vector $f^*$. Let $q\in\psi(Con_{n+1})$. Then $q$ satisfies as the inequality $(f,q)=(f^V,q^s)+(f^0,q^a)\ge 0$ (see (\ref{gF})) so the inequality $(f^*,q)=(f^V,q^s)-(f^0,q^a)\ge 0$. But it is easy to see that $(f,q)=(f^*,q^*)$ and $(f^*,q)=(f,q^*)$. This implies that $q^*\in\psi(Con_{n+1})$. \hfill $\Box$

\vspace{2mm}
The assertion of the following Proposition~\ref{Fas} is 
implied by the equality $(\psi(Con_{n+1}))^*=\psi(Con_{n+1})$.
 Call a facet $G$ of the cone $\psi(Con_{n+1})$ {\em symmetric}, if
 $q\in F$ implies $q^*\in F$, and call this facet  
{\em asymmetric}, if it is not symmetric.

\vspace{2mm}
\begin{prop}
\label{Fas}
Let $g\in Q_n$ be a facet vector of an asymmetric facet $G$ of the cone 
$\psi(Con_{n+1})$, and let $G^*=\{q^*:q\in G\}$. Then $G^*$ is a facet of 
$\psi(Con_{n+1})$, and $g^*$ is its facet vector.
\end{prop}

Recall that $Con_{n+1}$ has facets, that are zero-lifting of facets of $Con_n$. Call a facet $G$ of the cone $\psi(Con_{n+1})$ {\em zero-lifting} of a facet $F^V$ of $Con_n$ if $G=\psi(F)$, where $F$ is a facet of $Con_{n+1}$ which is zero-lifting of $F^V$.

\begin{prop}
\label{Fsy}
Let $g\in Q_n$ be a facet vector of a facet $G$ of the cone $\psi(Con_{n+1})$. Then the following assertions are equivalent:

(i) $g=g^*$;

(ii) the facet $G$ is symmetric;

(iii) $G=\psi(F)$, where $F$ is a facet of $Con_{n+1}$ which is zero-lifting of a facet $F^V$ of $Con_n$.

(iv) $G$ is a zero-lifting of a facet $F^V$ of $Con_n$.
\end{prop}
{\bf Proof}. (i)$\Rightarrow$(ii). If $g=g^*$, then $g=g^s$. Hence, $q\in G$ implies $(g,q)=(g^s,q)=(g^s,q^s)=(g,q^*)=0$. This means that $q^*\in G$, 
i.e., $G$ is symmetric.

(ii)$\Rightarrow$(i). By Proposition~\ref{iso}, the map $q\to q^*$ is an automorphism of $\psi(Con_{n+1})$. This map transforms a facet $G$ with facet vector $g$ into a facet $G^*$ with facet vector $g^*$. If $G$ is symmetric, then $G^*=G$, and therefore, $g^*=g$.

(iii)$\Rightarrow$(i). Let $f=f^V+f^0$ be a facet vector of a facet $F$ of $Con_{n+1}$ such that $f^0=0$. Then the facet $F$ is zero-lifting of the facet $F^V=F\cap{\mathbb R}^E$ of the cone $Con_n$. In this case, $f^V$ is also a facet vector of the facet $G=\psi(F)$ of $\psi(Con_{n+1})$. Obviously, $(f^V)^*=f^V$.

(iii)$\Rightarrow$(iv). This implication is implied by definition of zero-lifting of a facet of the cone $\psi(Con_{n+1})$.

(iv)$\Rightarrow$(i). The map $\chi$ induces a bijection between $\pi(F)$ and $\psi(F)$. Since $\pi(F)$ is zero-lifting of $F^V$, the facet vector of $\pi(F)$ belongs to ${\mathbb R}^E$. This implies that the facet vector $g$ of $\psi(F)$ belongs to ${\mathbb R}^E$, i.e., $g^*=g$. \hfill $\Box$

\vspace{2mm}
The symmetry group of $Con_{n+1}$ is the symmetric group $\Sigma_{n+1}$ of permutations of indices (see \cite{DL}). The group $\Sigma_n$ is a subgroup of the symmetry group of the cone $\psi(Con_{n+1})$. The full symmetry group of $\psi(Con_{n+1})$ is $\Sigma_n\times \Sigma_2$, where $\Sigma_2$ corresponds to the map $q\to q^*$ for $q\in\psi(Con_{n+1})$.  By Proposition~\ref{Fas}, the set of facets of $\psi(Con_{n+1})$ is partitioned into pairs $G,G^*$. But it turns out that there are pairs such that $G^*=\sigma(G)$, where $\sigma\in \Sigma_n$.

\section{Projections of hypermetric facets}
\label{prh}

The metric cone $Met_{n+1}$, the hypermetric cone $Hyp_{n+1}$ and the cut cone $Cut_{n+1}$ lying in the space ${\mathbb R}^{E\cup E_0}$ have an important class of {\em hypermetric} facets, that contains the class of {\em triangular} facets.

Let $b_i$, $i\in V$, be integers such that $\sum_{i\in V}b_i=\mu$, where $\mu=0$ or $\mu=1$. Usually these integers are denoted as a sequence $(b_1,b_2,...,b_n)$, where $b_i\ge b_{i+1}$. If, for some $i$, we have $b_i=b_{i+1}=...=b_{i+m-1}$, then the sequence is shortened as $(b_1,...,b_i^m,b_{i+m},...,b_n)$.

One relates to this sequence the following inequality of type $(b_1,...,b_n)$
\[(f(b),x)= -\sum_{i,j\in V}b_ib_jx_{(ij)}\ge 0, \]
where $x=\{x_{(ij)}\}\in{\mathbb R}^E$ and the vector $f(b)\in{\mathbb R}^E$ has coordinates $f(b)_{(ij)}=-b_ib_j$. This inequality is called of {\em negative} or {\em hypermetric} type if in the sum $\sum_{i\in V}b_i=\mu$ we have $\mu=0$ or $\mu=1$, respectively.

The set of hypermetric inequalities on the set $V\cup\{0\}$ determines a hypermetric cone $Hyp_{n+1}$. There are infinitely many hypermetric inequalities for metrics on $V\cup\{0\}$. But it is proved in \cite{DL}, that only finite number of these inequalities determines facets of $Hyp_{n+1}$. Since triangle inequalities are inequalities $(f(b),x)\ge 0$ of type $b=(1^2,0^{n-3},-1)$, the hypermetric cone $Hyp_{n+1}$ is contained in $Met_{n+1}$, i.e., $Hyp_{n+1}\subset Met_{n+1}$ with equality for $n=2$.

The hypermetric inequality $(f(b),x)\ge 0$ takes the following form on the set $V\cup\{0\}$.
\begin{equation}
\label{bV0}
-\sum_{i,j\in V\cup\{0\}}b_ib_jx_{(ij)}=-\sum_{(ij)\in E}b_ib_jx_{(ij)}-
\sum_{i\in V}b_0b_ix_{(0i)}\ge 0.
\end{equation}
If we decompose the vector $f(b)$ as $f(b)=f^V(b)+f^0(b)$, then $f^V(b)_{(ij)}=-b_ib_j$, $(ij)\in E$, and $f^0(b)_{(0i)}=-b_0b_i$, $i\in V$.

Let, for $S\subset V$, the equality $\sum_{i\in S}b_i=0$ hold. Denote by $b^S$ a sequence such that $b^S_i=-b_i$ if $i\in S$ and $b^S_i=b_i$ if $i\not\in S$. The sequence $b^S$ is called {\em switching} of $b$ by the set $S$.

The hypermetric cone $Hyp_{n+1}$ has the following property (see \cite{DL}). If an inequality $(f(b),x)\ge 0$ defines a facet and $\sum_{i\in S}b_i=0$ for some $S\subset V\cup\{0\}$, then the inequality $(f(b^S),x)\ge 0$ defines a facet, too.

{\bf Proof of Proposition~\ref{fso} for $Hyp_{n+1}$}.

Consider the inequality (\ref{bV0}), where $(f^0(b),e_0)=-\sum_{i\in V}b_0b_i=0$. Then $\sum_{i\in V}b_i=0$. Hence, the cone $Hyp_{n+1}$ has
similar inequality for $b^V$, where $b^V_i=-1$ for all $i\in V$. Hence, 
if one of these inequalities defines a facet, so does another. Obviously, $f^0(b^V)=-f^0(b)$. Hence, these facets satisfy the assertion of Proposition~\ref{fso}. \hfill $\Box$

\begin{theor}
\label{hyp}
Let $(f(b),x)\ge 0$ define a hypermetric facet of a cone in the space ${\mathbb R}^{E\cup E_0}$. Then the map $\psi$ transforms it either in a hypemetric facet if $b_0=0$ or in a distortion of a facet of negative type if $b_0=1$. Otherwise, the projection is not a facet.
\end{theor}
{\bf Proof}. By Section~\ref{psi}, the map $\psi$ transforms the hypermetric inequality (\ref{bV0}) for $x\in{\mathbb R}^{E\cup E_0}$ into the following inequality
\[-\sum_{(ij)\in E}b_ib_jq_{(ij)}-b_0\sum_{i\in V}b_iw_i\ge 0 \]
for $q=\sum_{(ij)\in E}q_{(ij)}\varphi(e_{(ij)})+\sum_{i\in V}w_iq(i)\in Q_n$.

Since $f(b)$ determines a hypermetric inequality, we have $b_0=1-\sum_{i\in V}b_i=1-\mu$. So, the above inequality takes the form
\[\sum_{(ij)\in E}b_ib_jq_{(ij)}\le (\mu-1)\sum_{i\in V}b_iw_i. \]

By Theorem~\ref{fac}, this facet is projected by the map $\psi$ into a facet if and only if $(f(b),e_0)=0$, where $e_0=\sum_{i\in V}e_{(0i)}$. Hence, we have
\[(f(b),e_0)=\sum_{i\in V}f(b)_{(0i)}=-\sum_{i\in V}b_0b_i=-b_0\mu=(\mu-1)\mu. \]
This implies that the hypermetric facet-defining inequality $(f(b),x)\ge 0$ is transformed into a facet-defining inequality if and only if either $\mu=0$ and then $b_0=1$ or $\mu=1$ and then $b_0=0$. So, we have

if $\mu=1$ and $b_0=0$, then the above inequality is a usual hypermetric inequality in the space $\psi({\mathbb R}^E)=\varphi({\mathbb R}^E)={\mathbb R}_s^{E^{\cal O}}$;

if $\mu=0$ and $b_0=1$, then the above inequality is the following distortion of an inequality of negative type
\begin{equation}
\label{neg}
-\sum_{(ij)\in E}b_ib_jq_{(ij)}-\sum_{i\in V}b_iw_i\ge 0, \mbox{  where  }\sum_{i\in V}b_i=0.
\end{equation}
\hfill $\Box$

\vspace{2mm}
Comparing (\ref{prod}) with the inequality (\ref{neg}), we see that a canonical facet vector $g(b)$ of a facet of $\psi(Hyp_{n+1})$ has the form $g(b)=g^s(b)+g^a(b)$, where $g_{ij}(b)=g_{(ij)}(b)+v_i-v_j$, and \[g_{(ij)}(b)=-b_ib_j, \mbox{   }v_i=-\frac{1}{n}b_i\mbox{  for all  }i\in V. \]

Define a cone of weighted quasi-hypermetrics $WQHyp_n=\psi(Hyp_{n+1})$
We can apply Proposition~\ref{iso}, in order to obtain the following assertion.
\begin{prop}
\label{hps}
The map $q\to q^*$ preserves the cone $WQHyp_n$, i.e.
\[(WQHyp_n)^*=WQHyp_n. \]
In other words, if $q\in WQHyp_n$ has weights $w_i, i\in V$, then the cone $WQHyp_n$ has a point $q^*$ with weights $-w_i, i\in V$. \hfill $\Box$
\end{prop}

\section{Generalizations of metrics}
The metric cone $Met_{n+1}$ is defined in the space ${\mathbb R}^{E\cup E_0}$. It has an extreme ray which is spanned by the vector $e_0=\sum_{i\in V}e_{(0i)}\in{\mathbb R}^{E_0}$. Facets of $Met_{n+1}$ are defined by the following set of triangle inequalities, where $d\in Met_{n+1}$.

Triangle inequalities of the sub-cone $Met_n$ that define facets of $Met_{n+1}$ that are zero-lifting and contain $e_0$:
\begin{equation}
\label{ijk}
d_{(ik)}+d_{(kj)}-d_{(ij)}\ge 0,\mbox{  for  }i,j,k\in V.
\end{equation}

Triangle inequalities defining facets that are not zero-lifting and contain the extreme ray $l_0$ spanned by the vector $e_0$:
\begin{equation}
\label{ij}
d_{(ij)}+d_{(j0)}-d_{(i0)}\ge 0\mbox{  and  }d_{(ij)}+d_{(i0)}-d_{(j0)}\ge 0,\mbox{  for  }i,j\in V.
\end{equation}

Triangle inequalities defining facets that do not contain the extreme ray $l_0$ and do not define facets of $Met_n$.
\begin{equation}
\label{ij0}
d_{(i0)}+d_{(j0)}-d_{(ij)}\ge 0,\mbox{  for  }i,j\in V.
\end{equation}

One can say that the cone $Met_n\in{\mathbb R}^{E_0}$ is lifted into
the space ${\mathbb R}^{E\cup E_0}$ using restrictions (\ref{ij}) and
(\ref{ij0}). Note that the inequalities (\ref{ij}) and (\ref{ij0}) imply the following inequalities of non-negativity
\begin{equation}
\label{pos}
 d_{(i0)}\ge 0,\mbox{  for  }i\in V.
\end{equation}

A cone defined by inequalities (\ref{ijk}) and (\ref{pos}) is called by
cone $WMet_n$ of {\em weighted} metrics $(d,w)$, where $d\in Met_n$ and $w_i=d_{(0i)}$ for $i\in V$ are {\em weights}.

If weights $w_i=d_{(0i)}$ satisfy additionally to the inequalities (\ref{pos}) also the inequalities (\ref{ij}), then the weighted metrics $(d,w)$ form a cone $dWMet_n$ of {\em down-weighted} metrics. If metrics have weights that satisfy the inequalities (\ref{pos}) and (\ref{ij0}), then these metrics are called {\em up-weighted} metrics. See details  in \cite{DD}, \cite{DDV}.

Above defined generalizations of metrics are functions on unordered pairs $(ij)\in E\cup E_0$. Generalizations of metrics as functions on ordered pairs $ij\in E^{\cal O}$ are called {\em quasi-metrics}.

The cone $QMet_n$ of quasi-metrics is defined in the space ${\mathbb R}^{E^{\cal O}}$ by non-negativity inequalities $q_{ij}\ge 0$ for all $ij\in E^{\cal O}$, and by triangle inequalities $q_{ij}+q_{jk}-q_{ik}\ge 0$ for all ordered triples $ijk$ for each $q\in QMet_n$. Below we consider in $QMet_n$ a sub-cone $WQMet_n$ of weighted quasi-metrics.

\section{Cone of weighted quasi-metrics}

We call a quasi-metric $q$ {\em weighted} if it belongs to the subspace $Q_n\subset{\mathbb R}^{E^{\cal O}}$. So, we define
\[ WQMet_n=QMet_n\cap Q_n. \]
A {\em quasi-metric} $q$ is called {\em weightable} if there are weights $w_i\ge 0$ for all $i\in V$ such that the following equalities hold
\[q_{ij}+w_i=q_{ji}+w_j \]
for all $i,j\in V$, $i\not=j$. Since $q_{ij}=q^s_{ij}+q^a_{ij}$, we have $q_{ij}+w_i=q^s_{ij}+q^a_{ij}+w_i=q^s_{ji}+q^a_{ji}+w_j$, i.e., $q^a_{ij}-q^a_{ji}=2q^a_{ij}=w_j-w_i$, what means that, up to multiple $\frac{1}{2}$ and sign, the antisymmetric part of $q_{ij}$ is $w_i-w_j$. So, weightable quasi-metrics are weighted.

Note that weights of a weighted quasi-metric are defined up to an additive constant. So, if we take weights non-positive, we obtain a weightable quasi-metric. Hence, sets of weightable and weighted quasi-metrics coincide.

By definition of the cone $WQMet_n$ and by symmetry of this cone, the triangle inequality $q_{ij}+q_{jk}-q_{ik}\ge 0$ and non-negativity inequality $q_{ij}\ge 0$ determine facets of the cone $WQMet_n$. Facet vectors of these facets are
\[t_{ijk}=e_{ij}+e_{jk}-e_{ik}\mbox{ and }e_{ij}, \]
respectively. It is not difficult to verify that $t_{ijk},e_{ij}\not\in Q_n$. 
Hence, these facet vectors are not canonical. Below, we give canonical representatives of these facet vectors.

Let $T(ijk)\subset K_n^{\cal O}$ be a triangle of $K_n^{\cal O}$ with direct arcs $ij,jk,ki$ and opposite arcs $ji,kj,ik$. Hence
\[f^{T(ijk)}=(e_{ij}+e_{jk}+e_{ki})-(e_{ji}+e_{kj}+e_{ik}). \]
\begin{prop}
\label{te}
Canonical representatives of facet vectors $t_{ijk}$ and $e_{ij}$ are
\[t_{ijk}+t^*_{ijk}=t_{ijk}+t_{kji},\mbox{  and  }g(ij)=(e_{ij}+e_{ji})+\frac{1}{n}(q(i)-q(j)), \]
respectively.
\end{prop}
{\bf Proof}. We have $t_{ijk}-f^{T(ijk)}=e_{ji}+e_{kj}-e_{ki}= t_{kji} =t^*_{ijk}$. This implies that the facet vectors $t_{ijk}$ and $t_{kji}$ determine the same facet, and the vector $t_{ijk}+t_{kji}\in{\mathbb R}_s^{E^{\cal O}}$ is a canonical representative of facet vectors of this facet. We obtain the first assertion of Proposition.

Consider now the facet vector $e_{ij}$. It is more convenient to 
take the doubled vector $2e_{ij}$. We show that the vector
\[g(ij)=2e_{ij}-\frac{1}{n}\sum_{k\in V-\{i,j\}}f^{T(ijk)} \]
is a canonical representative of the facet vector $2e_{ij}$. It is sufficient to show that $g(ij)\in Q_n$, i.e., $g_{kl}(ij)=g^s_{kl}(ij)+w_k-w_l$. In fact, we have $g_{ij}(ij)=2-\frac{n-2}{n}=1+\frac{2}{n}$, $g_{ji}(ij)=\frac{n-2}{n}=1-\frac{2}{n}$, $g_{ik}(ij)=-g_{ki}(ij)=\frac{1}{n}$, $g_{jk}(ij)=-g_{kj}(ij)=-\frac{1}{n}$, $g_{kk'}(ij)=0$. Hence, 
we have
\[g^s(ij)=e_{ij}+e_{ji}, \mbox{  }w_i=-w_j=\frac{1}{n},
\mbox{  and  }w_k=0\mbox{  for all  }k\in V-\{i,j\}. \]
These equalities imply the second assertion of Proposition. \hfill $\Box$

\vspace{2mm}
Let $\tau_{ijk}$ be a facet vector of a facet of $Met_n$ determined by the inequality $d_{(ij)}+d_{(jk)}-d_{(ik)}\ge 0$. Then $t_{ijk}+ t_{kji}=\varphi(\tau_{ijk})$, where the map $\varphi:{\mathbb R}^E\to {\mathbb R}_s^{E^{\cal O}}$ is defined in Section~\ref{spR}. Obviously, a triangular facet is symmetric.

Recall that $q_{ij}=q_{(ij)}+w_i-w_j$ if $q\in WQMet_n$. Let $i,j,k\in V$. It is not difficult to verify that the following equalities hold:
\begin{equation}
\label{met}
q^s_{ij}+q^s_{jk}-q^s_{ik}=q_{ij}+q_{jk}-q_{ij}\ge 0.
\end{equation}
Since $q^s_{ij}=q^s_{ji}=q_{(ij)}$, these inequalities show that the symmetric part $q^s$ of the vector $q\in WQMet_n$ is a semi-metric. Hence, if $w_i=w$ for all $i\in V$, then the quasi-semi-metric $q=q^s$ itself is a semi-metric. This implies that the cone $WQMet_n$ contains the cone of semi-metric $Met_n$. Moreover, $Met_n=WQMet_n\cap{\mathbb R}_s^{E^{\cal O}}$.

Now we show explicitly how the map $\psi$ transforms the cones 
$Met_{n+1}$ and $dWMet_n$ into the cone $WQMet_n$; see also Lemma 1 (ii) in \cite{DDV}.
\begin{theor}
\label{WQM}
The following equalities hold
\[\psi(Met_{n+1})=\psi(dWMet_n)=WQMet_n\mbox{  and  }WQMet^*_n=WQMet_n. \]
\end{theor}
{\bf Proof}. All facets of the metric cone $Met_{n+1}$ of metrics on the set $V\cup\{0\}$ are given by triangular inequalities $d_{(ij)}+d_{(ik)}-d_{(kj)}\ge 0$. They are hypermetric inequalities $(g(b),d)\ge 0$, where $b$ has only three non-zero values $b_j=b_k=1$ and $b_i=-1$ for some triple $\{ijk\}\subset V\cup\{0\}$. By Theorem~\ref{hyp}, the map $\psi$ transforms this facet into a hypermetric facet, i.e., into a triangular facets of the cone $\psi(Met_{n+1})$ if and only if $b_0=0$, i.e., if $0\not\in\{ijk\}$. If $0\in\{ijk\}$, then, by the same theorem, the equality $b_0=1$ should be satisfied. This implies $0\in\{jk\}$. In this case, the 
facet-defining inequality has the form (\ref{neg}), that in the case $k=0$, is
\[q_{(ij)}+w_i-w_j\ge 0. \]
This inequality is the non-negativity inequality $q_{ij}\ge 0$.

If $b_i=1,b_j=-1$ and $k=0$, the inequality $d_{(ij)}+d_{j0)}-d_{(0i)}\ge 0$ is transformed into inequality
\[q_{(ij)}+w_j-w_i\ge 0, \mbox{ i.e.,  }q^*_{ij}\ge 0. \]
This inequality and inequalities (\ref{met}) imply the last equality of this Theorem.

The inequalities (\ref{ij0}) define facets $F$ of $Met_{n+1}$ and 
$dWMet_n$ that do not contain the extreme ray $l_0$. Hence, by 
Theorem~\ref{hyp}, $\psi(F)$ are not facets of $WQMet_n$. But, 
recall that the cone $dWMet_n$ contains all facets of $Met_{n+1}$ excluding facets defined by the inequalities (\ref{ij0}). Instead of these facets, the cone $dWMet_n$ has facets $G_i$ defined by the non-negativity equalities (\ref{pos}) with facet vectors $e_{(0i)}$ for all $i\in V$. Obviously all these facets do not contain the extreme ray $l_0$. Hence, by Theorem~\ref{main}, $\psi(G_i)$ is not a facet of $\psi(dWMet_n)$. Hence, we have also the equality $WQMet_n=\psi(dWMet_n)$. \hfill $\Box$

\vspace{2mm}
{\bf Remark}. Facet vectors of facets of $Met_{n+1}$ that contain the extreme ray $l_0$ spanned by the vector $e_0$ are $\tau_{ijk}=\tau^V_{ijk}$, $\tau_{ij0}=\tau^V+\tau^0$ and $\tau_{ji0}=\tau^V-\tau^0$, where $\tau^V=e_{(ij)}$ and $\tau^0= e_{(j0)}-e_{(i0)}$. Hence, Proposition~\ref{fso} is true for $Met_{n+1}$, and we can apply Proposition~\ref{iso} in order to obtain the equality $WQMet_n^*=WQMet_n$ of Theorem~\ref{WQM}.

\section{The cone $Cut_{n+1}$}
The cut vectors $\delta(S)\in{\mathbb R}^{E\cup E_0}$ for all $S\subset V\cup\{0\}$ span all extreme rays of the cut cone $Cut_{n+1}\subset {\mathbb R}^{E\cup E_0}$. In other words, $Cut_{n+1}$ is the conic hull of all cut vectors. Since the cone $Cut_{n+1}$ is full-dimensional, its dimension is dimension of the space ${\mathbb R}^{E\cup E_0}$, that is $\frac{n(n+1)}{2}$.

Recall that $\delta(S)=\delta(V\cup\{0\}-S)$. Hence, we can consider only $S$ such that $S\subset V$, i.e., $0\not\in S$. Moreover, by Section~\ref{Cut},
\begin{equation}
\label{VS}
\delta(S)=\sum_{i\in S,j\not\in S}e_{(ij)}=\sum_{i\in S,j\in V-S}e_{(ij)}+\sum_{i\in S}e_{(0i)}=\delta^V(S)+\sum_{i\in S}e_{(0i)},
\end{equation}
where $\delta^V(S)$ is restriction of $\delta(S)$ on the space ${\mathbb R}^E=\psi({\mathbb R}^E)$. Note that
\[\delta(V)=\delta(\{0\})=\sum_{i\in V}e_{(0i)}=e_0. \]

Consider a facet $F$ of $Cut_{n+1}$. Let $f$ be facet vector of $F$. Set
\[R(F)=\{S\subset V:(f,\delta(S))=0\}.\]
For $S\in R(F)$, the vector $\delta(S)$ is called {\em root} of the facet $F$. By (\ref{VS}), for $S\in R(F)$, we have
\begin{equation}
\label{fS}
(f,\delta(S))=(f,\delta^V(S))+\sum_{i\in S}f_{(0i)}=0.
\end{equation}
We represent each facet vector of $Cut_{n+1}$ as $f=f^V+f^0$, where $f^V\in{\mathbb R}^E$ and $f^0\in{\mathbb R}^{E_0}$.

The set of facets of the cone $Cut_{n+1}$ is partitioned onto equivalence classes by {\em switchings} (see \cite{DL}). For each $S,T\subset V\cup\{0\}$, the switching by the set $T$ transforms the cut vector $\delta(S)$ into the vector $\delta(S\Delta T)$, where $\Delta$ is symmetric difference, i.e., $S\Delta T=S\cup T-S\cap T$. It is proved in \cite{DL} that if $T\in R(F)$, then $\{\delta(S\Delta T):S\in R(F)\}$ is the set of roots of the switched facet $F^{\delta(T)}$ of $Cut_{n+1}$. Hence, $R(F^{\delta(T)})=\{S\Delta T:S\in R(F)\}$.

Let $F$ be a facet of $Cut_{n+1}$. Then $F$ contains the vector $e_0=\delta(V)$ if and only if $V\in R(F)$. Hence, Lemma~\ref{FV} below is an extended reformulation of Proposition~\ref{fso}.
\begin{lem}
\label{FV}
Let $F$ be a facet of $Cut_{n+1}$ such that $V\in R(F)$. Let $f=f^V+f^0$ be facet vector of $F$. Then the vector $f^*=f^V-f^0$ is facet vector of
switching $F^{\delta(V)}$ of the facet $F$, and $V\in R(F^{\delta(V)})$.
\end{lem}
{\bf Proof}. Since $V\in R(F)$, $F^{\delta(V)}$ is a facet of $Cut_{n+1}$. Since $S\Delta V=V-S={\overline S}$, for $S\subset V$, we have
\[R(F^{\delta(V)})=\{{\overline S}:S\in R(F)\}. \]
Since $\emptyset\in R(F)$, the set $\emptyset\Delta V=V\in R(F^{\delta(V)})$. Now, using (\ref{fS}), for $S\in R(F^{\delta(V)})$, we have
\[(f^*,\delta(S))=((f^V-f^0),\delta(S))=(f^V,\delta^V(S))-\sum_{i\in S}f_{(0i)}.\]
Note that $\delta^V({\overline S})=\delta^V(S)$, and, since $V\in R(F)$, $\delta(V)=\delta(\{0\})$, we have $(f,\delta(V))=\sum_{i\in V}f_{(0i)}=0$. 
Hence, $\sum_{i\in{\overline S}}f_{(0i)}=-\sum_{i\in S}f_{(0i)}$. It is easy to see, that $(f^*,\delta(S))=(f,\delta({\overline S}))$. Since $S\in R(F^{\delta(V)})$ if and only if ${\overline S}\in R(F)$, we see that $f^*$ is a facet vector of $F^{\delta(V)}$. \hfill $\Box$

\vspace{2mm}
The set of facets of $Cut_{n+1}$ is partitioned into orbits under action of the permutation group $\Sigma_{n+1}$. But some permutation non-equivalent facets are equivalent under switchings. We say that two facets $F,F'$ of $Cut_{n+1}$ belong to the same {\em type} if there are $\sigma\in \Sigma_{n+1}$ and $T\subset V$ such that $\sigma(F')=F^{\delta(T)}$.

\section{Cone $OCut_n$}
Denote by $OCut_n\subset{\mathbb R}^{E^{\cal O}}$ the cone whose extreme rays are spanned by ocut vectors $c(S)$ for all $S\subset V$, $S\not=\emptyset, V$. In other words, let
\[OCut_n=\{c\in Q_n:c=\sum_{S\subset V}\alpha_Sc(S), \mbox{  }\alpha_S\ge 0\}. \]
Coordinates $c_{ij}$ of a vector $c\in OCut_n$ are given in (\ref{cw}), where $w_i\ge 0$ for all $i\in V$. Hence, $OCut_n\subset Q_n$. Recall that
\begin{equation}
\label{cqS}
c(S)=\frac{1}{2}(\delta^{\cal O}(S)+\sum_{i\in S}q(i)),
\end{equation}
where $\delta^{\cal O}(S)=\varphi(\delta^V(S))$. Note that $\delta^{\cal O}({\overline S})=\delta^{\cal O}(S)$ and $q({\overline S})=-q(S)$, where ${\overline S}=V-S$.

Denote by $Cut_n^{\cal O}=\varphi(Cut_n)$ the cone generated by $\delta^{\cal O}(S)$ for all $S\subset V$. The vectors $\delta^{\cal O}(S)$ for all $S\subset V$, $S\not=\emptyset, V$, are all extreme rays of the cone $Cut_n^{\cal O}$ that we identify with $Cut_n$ embedded into the space ${\mathbb R}^{E^{\cal O}}$.

\begin{lem}
\label{cS}
For $S\subset V$, the following equality holds:
\[\psi(\delta(S))=2c(S). \]
\end{lem}
{\bf Proof}. According to Section~\ref{psi}, $\psi(\delta^V(S))=\varphi(\delta^V(S))=\delta^{\cal O}(S)$. Besides, $\psi(e_{(0i)})=q(i)$ for all $i\in V$. Hence, using (\ref{VS}), we obtain
\[\psi(\delta(S))=\psi(\delta^V(S))+\sum_{i\in S}\psi(e_{(0i)})=\varphi(\delta^V(S))+\sum_{i\in S}q(i)=\delta^{\cal O}(S)+q(S). \]

Recall that $\psi(\delta(V))=\psi(e_0)={\bf 0}$ and $c(V)=0$. Hence, according to (\ref{cqS}), we obtain
\[\psi(\delta(S))=2c(S),\mbox{  for all  }S\subset V.\]
Lemma is proved. \hfill $\Box$

\begin{theor}
\label{psC}
The following equalities hold:
\[\psi(Cut_{n+1})=OCut_n\mbox{  and  }OCut^*_n=OCut_n.\]
\end{theor}
{\bf Proof}. Recall that the conic hull of vectors $\delta(S)$ for all $S\subset V$ is $Cut_{n+1}$. The conic hull of vectors $c(S)$ for all $S\subset V$ is the cone $OCut_n$. Since $\psi(\delta(V))=c(V)={\bf 0}$, the first result follows.

The equality $OCut^*_n=OCut_n$ is implied by the equalities $c^*(S)=c({\overline S})$ for all $S\subset V$.

By Lemma~\ref{FV}, the equality $OCut^*_n=OCut_n$ is the special case $Con_{n+1}=Cut_{n+1}$ of 
Proposition~\ref{iso}.
\hfill $\Box$

\section{Facets of $OCut_n$}

\begin{lem}
\label{fV}
Let $F$ be a facet of $Cut_{n+1}$. Then $\psi(F)$ is a facet of $OCut_n$ if and only if $V\in R(F)$.
\end{lem}
{\bf Proof}. By Theorem~\ref{main}, $\psi(F)$ is a facet of $OCut_n$ if and only if $e_0=\delta(V)\subset F$, i.e., if and only if $V\in R(F)$. \hfill $\Box$

\vspace{2mm}
For a facet $G$ of $OCut_n$ with facet vector $g$, we set
\[R(G)=\{S\subset V:(g,c(S))=0\}\]
and call the vector $c(S)$ for $S\in R(G)$ by {\em root} of the facet $G$.

Note that $\delta(\emptyset)={\bf 0}$ and $c(\emptyset)=c(V)={\bf 0}$. Hence, $\emptyset\in R(F)$ and $\emptyset\in R(G)$ for all facet $F$ of $Cut_{n+1}$ and all facets $G$ of $OCut_n$. The roots $\delta(\emptyset)={\bf 0}$ and $c(\emptyset)=c(V)={\bf 0}$ are called {\em trivial} roots.
\begin{prop}
\label{RG}
For a facet $F$ of $Cut_{n+1}$, let $G=\psi(F)$ be a facet of $OCut_n$. 
Then the following equality holds:
\[R(G)=R(F).\]
\end{prop}
{\bf Remark}. We give two proofs of this equality. Both are useful.

{\bf First proof}. According to Section~\ref{psi}, the map $\psi$ transforms an inequality $(f,x)\ge 0$ defining a facet of $Cut_{n+1}$ into the inequality (\ref{gF}) defining the facet $G=\psi(F)$ of $OCut_n$. Recall the the inequality (\ref{gF}) relates to the representation of vectors $q\in Q_n$ in the basis $\{\varphi(e_{ij}),q(i)\}$, i.e., $q=\sum_{(ij)\in E}q_{(ij)}\varphi(e_{(ij)})+\sum_{i\in V}w_iq(i)$. Let $q=c(S)$ for $S\in R(G)$. Then, according to (\ref{cqS}), we have $q_{(ij)}=\frac{1}{2}\delta^V_{(ij)}(S)$, $w_i=\frac{1}{2}$ for $i\in S$ and $w_i=0$ for $i\in{\overline S}$. Hence, omitting the multiple $\frac{1}{2}$, the inequality in (\ref{gF}) gives the following equality
\[\sum_{(ij)\in E}f_{(ij)}\delta^V_{(ij)}(S)+\sum_{i\in S}f_{(0i)}=0 \]
which coincides with (\ref{fS}). This implies the assertion of this Proposition.

{\bf Second proof}. By Theorem~\ref{main}, $\psi(l)$ is an extreme ray of $\psi(F)$ if and only if $l$ is an extreme ray of $F$ and $l\not=l_0$. Since $l$ is spanned by $\delta(S)$ for some $S\in R(F)$ and $\psi(l)$ is spanned by $\psi(\delta(S))=c(S)$, we have $R(G)=\{S\subset V:S\in R(F)\}$. Since $c(V)={\bf 0}$, we can suppose that $V\in R(G)$, and then $R(G)=R(F)$. \hfill $\Box$

\vspace{2mm}
{\bf Remark}. Note that $\delta(V)=\delta(\{0\})=e_0\not={\bf 0}$ is a non-trivial root of $F$, i.e., $V\in R(F)$. But $c(V)=\psi(\delta(V))={\bf 0}$ is a trivial root of $R(G)$.

\vspace{2mm}
Recall that, for a subset $T\subset V$, we set ${\overline T}=V-T$. Note that ${\overline T}=V\Delta T$ and ${\overline T}\not=V\cup\{0\}-T$.
\begin{lem}
\label{FT}
Let $F$ be a facet of $Cut_{n+1}$, and $T\in R(F)$. Then the image $\psi(F^{\delta(T)})$ of the switched facet $F^{\delta(T)}$ is a facet of $OCut_n$ if and only if ${\overline T}\in R(F)$.
\end{lem}
{\bf Proof}. By Lemma~\ref{fV}, $\psi(F^{\delta(T)})$ is a facet of $OCut_n$ if and only if $V\in R(F^{\delta(T)})$, i.e., if and only if $V\Delta T={\overline T}\in R(F)$. \hfill $\Box$

\vspace{2mm}
For a facet $G$ of $OCut_n$, define $G^{\delta(T)}$ as the conic hull of $c(S\Delta T)$ for all $S\in R(G)$. Since each facet $G$ of $OCut_n$ is $\psi(F)$ for some facet $F$ of $Cut_{n+1}$, Lemma~\ref{FT} and Proposition~\ref{RG} imply the following assertion.
\begin{theor}
\label{swi}
Let $G$ be  a facet of $OCut_n$. Then $G^{\delta(T)}$ is a facet of $OCut_n$ if and only if $T,{\overline T}\in R(G)$, and then $R(G^{\delta(T)})=\{S\Delta T:S\in R(G)\}$. \hfill $\Box$
\end{theor}
Theorem~\ref{swi} asserts that the set of facets of the cone $OCut_n$ is partitioned onto equivalence classes by switchings $G\to G^{\delta(T)}$, where $T,{\overline T}\in R(G)$.

The case $T=V$ in Theorem~\ref{swi} plays a special role.
Recall that $V\in R(F)$ if $F$ is a facet of $Cut_{n+1}$ such that $\psi(F)$ is a facet of $OCut_n$. Hence, Lemma~\ref{FV} and Proposition~\ref{iso} imply the following fact.
\begin{prop}
\label{gf}
Let $F$ be a facet of $Cut_{n+1}$ such that $\psi(F)$ is a facet of $OCut_n$. Let $g=g^s+g^a$ be a facet vector of the facet $\psi(F)$. Then the vector $g^*=g^s-g^a$ is a facet vector of the facet $\psi(F^{\delta(V)})=(\psi(F))^*=(\psi(F))^{\delta(V)}$ such that $R((\psi(F))^*)=\{{\overline S}:S\in R(F)\}$. \hfill $\Box$
\end{prop}

Recall that roughly speaking $OCut_n$ is projection of $Cut_{n+1}$ along the vector $\delta(V)=\delta(\{0\})$.

Let $\sigma\in\Sigma_n$ be a permutation of the set $V$. For a vector $q\in {\mathbb R}^{E^{\cal O}}$, we have  $\sigma(q)_{ij}=q_{\sigma(i)\sigma(j)}$. Obviously if $g$ is a facet vector of a facet $G$ of $OCut_n$, then $\sigma(g)$ is the facet vector of the facet $\sigma(G)=\{\sigma(q):q\in G\}$.

Note that, by Proposition~\ref{gf}, the switching by $V$ is equivalent to the operation $q\to q^*$. Hence, the symmetry group of $OCut_n$ contains the group $\Sigma_n\times\Sigma_2$, where $\Sigma_2$ relates to the map $q\to q^*$ for $q\in OCut_n$.
\begin{theor}
\label{sym}
The group $\Sigma_n\times\Sigma_2$ is the symmetry group of the cone $OCut_n$.
\end{theor}
{\bf Proof}. Let $\gamma$ be a symmetry of $OCut_n$. Then $\gamma$ is a symmetry of the set ${\cal F}(e_0)$ of facets $F$ of the cone $Cut_{n+1}$ containing the vector $e_0$. The symmetry group $\Gamma(e_0)$ of the set ${\cal F}(e_0)$ is a subgroup of the symmetry group of the cut-polytope $Cut_{n+1}^{\Box}$. In fact, $\Gamma(e_0)$ is stabilizer of the edge $e_0$ of the polytope $Cut_{n+1}^{\Box}$. But it is well-known that $\Gamma(e_0)$ consists of the switching by $V$ and permutations $\sigma\in\Sigma_{n+1}$ leaving the edge $e_0$ non-changed. The map $\psi$ transforms these symmetries of ${\cal F}(e_0)$ into symmetries $\sigma\in\Sigma_n$ and $q\to q^*$ of the cone $OCut_n$. \hfill $\Box$

\vspace{2mm}
The set of all facets of $OCut_n$ is partitioned onto orbits of facets that are equivalent by the symmetry group $\Sigma_n\times\Sigma_2$. It turns out that, for some facets $G$, subsets $S\in R(G)$ and permutations $\sigma\in\Sigma_n$, we have $G^{\delta(S)}=\sigma(G)$.

By Proposition~\ref{Fsy}, if a facet of $Cut_{n+1}$ is zero-lifting of a facet $F^V$ of $Cut_n$, then the facet $G=\psi(F)$ of $OCut_n$ is symmetric and $G=G^*=G^{\delta(V)}$ is zero-lifting of $F^V$.

So, there are two important classes of orbits of facets of $OCut_n$. Namely, the orbits of symmetric facets, that are zero-lifting of facets of $Cut_n$, and orbits of asymmetric facets that are $\psi$-images of facets of $Cut_{n+1}$ and are not zero-lifting.

\section{Cases $3\le n\le 6$} 
Compare results of this Section with Table 2 of \cite{DDV}. 

Most of described below facets are hypermetric or negative type. We give here the corresponding vectors $b$ in accordance with Section~\ref{prh}. 

{\bf n=3}. Note that $Cut_4=Hyp_4=Met_4$. Hence,
\[OCut_3=WQHyp_3=WQMet_3. \]
All these cones have two orbits of facets: one orbit of non-negativity facets with $b=(1,0,-1)$ and another orbit of triangular facets with $b=(1^2,-1)$.

{\bf n=4}. We have $Cut_5=Hyp_5\subset Met_5$. Hence,
\[OCut_4=WQHyp_4\subset WQMet_4.\]
The cones $Hyp_5=Cut_5$ have two orbits of facets: triangular and pentagonal facets. Recall that a triangular facet with facet vector $\tau_{ijk}$ is zero-lifting if $0\not\in\{ijk\}$. Hence, the cones $WQHyp_4=OCut_4$ have three orbits of facets: of non-negativity with $b=(1,0^2,-1)$, triangular with $b=(1^2,0,-1)$ and weighted version of negative type with $b=(1^2,-1^2)$.

{\bf n=5}. We have again $Cut_6=Hyp_6\subset Met_6$. Hence,
\[OCut_5=WQHyp_5\subset WQMet_5.\]
The cones $Hyp_6=Cut_6$ have four orbits of facets, all are hypermetric: triangular with $b=(1^2,0^3,-1)$, pentagonal with $b=(1^3,0,-1^2)$ and two more types, one with $b=(2,1^2,-1^3)$ and its switching with $b=(1^4,-1,-2)$. These four types provide 6 orbits of facets of the cones $WQHyp_5=OCut_5$: non-negativity with $b=(1,0^3,-1)$, triangular with $b=(1^2,0^2,-1)$, of negative type with $b=(1^2,0,-1^2)$, pentagonal with $b=(1^3,-1^2)$, and two of negative type with $b=(2,1,-1^3)$ and $b=(1^3,-1,-2)$. 

The last two types belong to the same orbit of the full symmetry group $\Sigma_5\times\Sigma_2$. Hence, the cone $OCut_5$ has 5 orbits of facets under action of its symmetry group. 

{\bf n=6}. Now, we have $Cut_7\subset Hyp_7\subset Met_7$. Hence
\[OCut_6\subset WQHyp_6\subset WQMet_6.\]

The cone $Cut_7$ has 36 orbits of facets that are equivalent under action of the permutation group $\Sigma_7$. Switchings contract these orbits into 11 types $F_k$, $1\le k\le 11$, of facets that are switching equivalent (see \cite{DL}, Sect. 30.6). J.Vidali compute orbits of facets of $OCut_6$ under action of the group $\Sigma_6$. Using these computations, we give in Table below numbers of orbits of facets of cones $Cut_7$ and $OCut_6$ (cf. Figure 30.6.1 of \cite{DL}).

The first row of Table gives types of facets of $Cut_7$. In the second row of Table, for each type $F_k$, numbers of orbits of facets of $Cut_7$ of type $F_k$ under action of the group $\Sigma_7$. The third row of Table, for each type $F_k$, gives numbers of orbits of facets of $OCut_6$ that are obtained from facets of type $F_k$ under action of the group $\Sigma_6$. The fourth row gives, for each type $F_k$, numbers of orbits of facets of $OCut_6$ that are obtained from facets of type $F_k$ under action of the group $\Sigma_6\times\Sigma_2$.

The last column of Table gives total numbers of orbits of facets of the cones $Cut_7$ and $OCut_6$.

{\bf Table}.
\[\begin{array}{|l||l|l|l||l|l|l||l|l|l||l|l||l|} \hline
\mbox{types} &F_1 &F_2 &F_3 &F_4 &F_5 &F_6 &F_7 &F_8 &F_9 &F_{10} &F_{11}&|\Omega|\\ \hline
\Sigma_7 &1 &1 &2 &1 &3 &2 &4 &7 &5 &3 &7 & 36\\ \hline \hline
\Sigma_6 &2 &2 &4 &1 &3 &2 &7 &13 &6 &6 &15 &61\\ \hline
\Sigma_6\times\Sigma_2 &2 &2 &3 &1 &2 &1 &4 &7 &3 &4 &8 &37\\ \hline
\end{array} \]

\vspace{2mm}

The first three types $F_1,F_2,F_3$ relate to 4 orbits of hypermetric facets $F(b)$ of $Cut_7$ that are zero-lifting, where $b=(1^2,0^4,-1)$, $b=(1^3,0^2,-1^2)$ and $b=(2,1^2,0,-1^3)$, $b=(1^4,0,-1,-2)$. Each of these four orbits of facets of $Cut_7$ under action of $\Sigma_7$ gives two orbits of facets of $OCut_6$ under action of the group $\Sigma_6$.

The second three types $F_4,F_5,F_6$ relate to 6 orbits of hypermetric facets $F(b)$ of $Cut_7$ that are not zero-lifting. Each of these 6 orbits gives one orbit of facets of $OCut_6$ under action of the group $\Sigma_6$.

The third three types $F_7,F_8,F_9$ relate to 16 orbits of facets of clique-web types $CW^7_1(b)$. These 16 orbits give 26 orbits of facets of $OCut_6$ under action of $\Sigma_6$. 

The last two types $F_{10}=Par_7$ and $Gr_7$ are special (see \cite{DL}). They relate to 10 orbits of $Cut_7$, that give 21 orbits of facets of $OCut_6$ under action of $\Sigma_6$.

The subgroup $\Sigma_2$ of the full symmetry group $\Sigma_6\times\Sigma_2$ contracts some pairs of orbits of the group $\Sigma_6$ into one orbit of the full group. The result is given in the forth row of Table. 

Note that the symmetry groups of $Cut_7$ and $OCut_6$ have 36 and 37 orbits of facets, respectively.


\newpage
\begin{thebibliography}{99}

\bibitem[AACMP97]{AAC}
O.Aichholzer, F.Auerenhammer, D.Z.Chen, D.T.Lee, A.Mukhopadhyay and E.Papadopoulou, {\em Voronoi diagrams for direction-sensitive distances}, Proceedings of $13^{th}$ Annual ACM Simposium Computational Geometry, Nice, France (1997), 418--420.

\bibitem[Aig79]{Aig}
M.Aigner, {\em Combinatorial Theory}, Springer-Verlag, Berlin 1979, (Chapter VII 3B).

\bibitem[AM11]{AM}
D.Avis and C.Meagher, {\em On the Directed Cut Cone and Polytope}, Manuscript Draft No. MAPR-D-11-00057, 2011.

\bibitem[Bir67]{Bir}
G.Birkhoff, {\em Lattice Theory}, AMS, Providence, Rhode Island, 1967.

\bibitem[CMM06]{CMM}
M.Charikar, K.Makarychev, and Y.Makarychev, {\em Directed metrics and directed graph partitioning problem} Proc. of 17th ACM-SIAM Symposium on Discrete Algorithms (2006) 51--60.

\bibitem[DD10]{DD}
M.Deza and E.Deza, {\em Cones of Partial  Metrics}, Contributions in Discrete Mathematics, {\bf 6} (2010), 26--41.

\bibitem[DDV11]{DDV}
M.Deza, E.Deza and J.Vidali, {\em Cones of Weighted and Partial Metrics}, 
arXiv: 1101.0517v2[math.Co]04 Jan  2011, to appear in {\em Algebra 2010: Advances in Algebraic Structures}, World Scientific, 2011.

\bibitem[DL97]{DL}
M.Deza and M.Laurent, {\em Geometry of cuts and metrics}, Springer-Verlag, 1997.

\bibitem[Ha14]{Ha}
F. Hausdorff, {\em Grundz\"uge der Mengenlehre}, Leipzig, Verlag "Veit and Co", 1914.

\bibitem[Se97]{Se}
A.K.Seda, {\em Quasi-metrics and semantic of logic programs}, Fundamenta Informaticae, {\bf 29} (1997), 97--117.

\end{thebibliography}
\end{document}